\documentclass[12pt,draft,a4paper]{amsart}
\usepackage{amssymb}
\usepackage{amsmath}
\usepackage{amsfonts}

\newcommand{\C}{\mathbb{C}}

\DeclareMathOperator{\Irr}{Irr}
\newtheorem{main-dummy}{Main-Dummy}
\newtheorem{dummy}{Dummy}

\newtheorem{main-theorem}[main-dummy]{Theorem}

\numberwithin{dummy}{section}
\numberwithin{equation}{section}

\newtheorem{theorem}[dummy]{Theorem}

\newtheorem{cor}[dummy]{Corollary}

\theoremstyle{definition}

\theoremstyle{remark}
\newtheorem{rem}[dummy]{Remark}
\begin{document}
\bibliographystyle{amsalpha}

%%%%%%%%%%%%%%%%%%%%%%%%%%%%%%%%%%%%%%%%%%%%%%%%%%%%%%%%%%%%%%%%%%%%%%%
%%%%%%%%%%%%%%%%%%%%%%%%%%%%%%%% Title %%%%%%%%%%%%%%%%%%%%%%%%%%%%%%%%
\author{Sandro Mattarei}

\email{mattarei@science.unitn.it}

\urladdr{http://www-math.science.unitn.it/\~{ }mattarei/}

\address{Dipartimento di Matematica\\
  Universit\`a degli Studi di Trento\\
  via Sommarive 14\\
  I-38050 Povo (Trento)\\
  Italy}

\title[Retrieving information about a group]%
{Retrieving information about a group from its character degrees or from its class sizes}

\begin{abstract}
We prove that a knowledge of the character degrees of a finite group $G$ and of their multiplicities
determines whether $G$ has a Sylow $p$-subgroup as a direct factor.
An analogous result based on a knowledge of the conjugacy class sizes was known.
We prove variations of both results and discuss their similarities.
\end{abstract}

\date{31 March 2005}



\thanks{The  author  is grateful  to  Ministero dell'Istruzione, dell'Universit\`a  e
  della  Ricerca, Italy,  for  financial  support of the
  project ``Graded Lie algebras  and pro-$p$-groups of finite width''.}

\maketitle

\thispagestyle{empty}
%%%%%%%%%%%%%%%%%%%%%%%%%%%%%%%%%%%%%%%%%%%%%%%%%%%%%%%%%%%%%%%%%%%%%%%
%%%%%%%%%%%%%%%%%%%%%%%%%%%%%%%%%%%%%%%%%%%%%%%%%%%%%%%%%%%%%%%%%%%%%%%

%----------------------------------------------------------------
\section{Introduction}\label{sec:intro}

Several results in the literature assert that certain properties of a finite group $G$ can or cannot
be inferred from knowledge of the complex group algebra $\C G$ (as a $\C$-algebra).
This piece of information is equivalent to knowing the sizes of the full matrix algebras into which
$\C G$ decomposes according to Wedderburn's theorem or, equivalently,
the {\em character degree frequency function} $m_G$ defined as
\[
m_G(n)=|\{\chi\in\Irr(G):\chi(1)=n\}|.
\]
For example, Isaacs showed in~\cite{Isaacs:group_algebra} that if $\C G\cong\C H$ for finite groups $G$ and $H$
then the nilpotency of $G$ implies the nilpotency of $H$.
Hawkes exhibited examples in~\cite{Hawkes-supersoluble} which show that this is not the case with supersolubility
replacing nilpotency.
The positive results usually give an explicit criterion for checking the validity of the desired property
in terms of $m_G$.
In particular, we formalize as follows an invariant of $\C G$ derived from $m_G$ which seems to often play a role.
For an integer $n$ and a set $\pi$ of primes let $n=n_{\pi}\cdot n_{\pi'}$ be its decomposition into the product
of its $\pi$- and $\pi'$-parts.
If $G$ is a finite group we set
\[
u_\pi(G)=
\sum_{\text{$n_{\pi'}=1$}} m_G(n)\cdot n^2=
\sum_{\substack{\chi\in\Irr(G)\\ \text{$\chi(1)_{\pi'}=1$}}}
\chi(1)^2.
\]

For example, it is known that $G$ has a normal abelian Hall $\pi$-subgroup
if and only if all character degrees of $G$ are $\pi'$-numbers, that is,
$u_{\pi'}(G)=|G|$
or, equivalently,
$u_{\pi}(G)=|G:G'|$.
The necessity of the condition follows from Ito's Theorem~\cite[(6.15)]{Isa}.
The sufficiency follows easily from the special case $\pi=\{p\}$, which is
also due to Ito for $G$ soluble~\cite[(12.33)]{Isa}
and has been extended to arbitrary finite groups by G.~Michler
using the classification of finite simple groups (see~\cite[Theorem~5.4]{Michler:defects}).

The result of Isaacs mentioned above on how to read nilpotency of $G$ from $\C G$ involves this invariant
in a less trivial way.
In fact, Isaacs proved in~\cite[Theorem~1]{Isaacs:group_algebra} that
$|G:G'|_p$ divides $u_{p'}(G)$ for any finite group $G$, and that $G$ is $p$-nilpotent
(that is, $G$ has a normal $p$-complement)
if and only if the quotient of these two numbers is not divisible by $p$.
In particular, whether $G$ is nilpotent can be read from $m_G$,
because $G$ is nilpotent if and only if it is $p$-nilpotent for every prime $p$ dividing $|G|$.
In short, Isaacs' condition can be written as
$u_{p'}(G)_p=|G:G'|_p$.

It seems to be an open question whether similar results may hold with $\pi$ (or $\pi'$, if we like)
replacing $p$.
In particular, it is not known whether the existence of a normal Sylow $p$-subgroup of $G$
can be read from the function $m_G$.
One of the goals of this note is to prove that the question whether $G$ has a normal Sylow $p$-subgroup
{\em and} a normal $p$-complement, that is, whether $G$ is the direct product of a $p$-subgroup
and a $p'$-subgroup, can be answered in terms of $m_G$.
We prove in Theorem~\ref{thm:direct-product-Isaacs} that this occurs if and only if $u_{p'}(G)=|G|_{p'}\cdot|G:G'|_{p}$.

Another equivalent condition for $p$-nilpotency was found by Cossey and Hawkes,
but in terms of the character degrees of $G$ which are a power of $p$.
They proved in~\cite[Theorem~1]{Cossey-Hawkes:group_algebra} that $G$ is $p$-nilpotent if and only if
$u_{p}(G)_p=|G|_p$.
Actually, Cossey and Hawkes showed, more generally, that
$u_{p}(G)_p=|G:\mathbf{O}^p(G)|$ for every finite group $G$.
Consequently,
the index $|G:K_\infty(G)|$ of the last term of the lower central series of $G$
equals the product $\prod_p u_{p}(G)_p$ over all primes $p$
(or just the prime divisors of $|G:G'|$, since $K_\infty(G)=\bigcap_{p\mid|G:G'|}\mathbf{O}^p(G)$).
Since this more general result is apparently not attainable by Isaacs' method,
one may get the impression that $u_{p}(G)$ contains more precise information than $u_{p'}(G)$
concerning the $p$-structure of $G$.
This impression is somehow in contrast with our next result, Theorem~\ref{thm:direct-product-Cossey-Hawkes}, where
$u_{p}(G)=|G|_{p}\cdot|G:G'|_{p'}$
is shown to occur if and only if
$G$ has a normal $p$-complement $N$ with $[N,G]=N'$.
In order to conclude that $G$ is the direct product of a $p$-group and a $p'$-group
it is then sufficient to have any additional information which forces $N$ to be nilpotent,
for instance as in Corollary~\ref{cor:direct-product-Cossey-Hawkes}.

Several results which allow one to deduce properties of a finite group from its
character degrees have
counterparts with the conjugacy class sizes replacing the character degrees.
It is natural to expect that class sizes give information on the group somehow ``from below''
({\em e.g.}~about subgroups),
being kind of dual to the character degrees, which give information ``from above''
({\em e.g.}~about quotient groups).
The results based on class sizes are usually easier to prove than their analogues with character degrees,
and often correspond to stronger structural conditions on the group.
A perhaps subjective interpretation of this difference is that, rather than to the character degrees,
the class sizes are dual to the squares of the character degrees.
The latter have the same average value as the class sizes,
but being squares they are bound to have less arithmetic content.

For example, it is known (and very easy to prove, see~\cite[Theorem~5']{Huppert:Mainz}) that
$G$ has a central Sylow $p$-subgroup if and only if
all class sizes are $p'$-numbers.
This generalizes at once to $\pi$ instead of $p$, and is the class-size analogue of the Ito-Michler theorem mentioned above about
recognizing abelian Hall $\pi$-subgroups from the character degrees.

Less trivially, Cossey, Hawkes and Mann proved in~\cite{Cossey-Hawkes-Mann} that the order of the hypercentre of $G$
can be computed from the {\em class size frequency function}
\[
w_G(n)=\frac{1}{n}|\{g\in G: |G:\mathbf{C}_G(g)=n|\}|.
\]
More precisely, they proved that $|Z_\infty(G)|_p=|\mathcal{S}_p(G)|_p$, where $\mathcal{S}_p(G)$ is the union of those classes
of $G$ whose cardinality is a power of $p$.
In particular, whether $G$ is nilpotent can be read from $w_G$.
We generalize their notation to a set of primes $\pi$ and denote by $\mathcal{S}_\pi(G)$ the union of those classes
of $G$ whose cardinality is a $\pi$-number.
Thus,
\[
|\mathcal{S}_\pi(G)|=\sum_{\text{$n_{\pi'}=1$}} w_G(n)\cdot n
\]
is an analogue for class sizes of the function $u_\pi(G)$ for character degrees,
and the result of Cossey, Hawkes and Mann~\cite{Cossey-Hawkes-Mann} considered here
is the dual version of the result of Cossey and Hawkes~\cite[Theorem~1]{Cossey-Hawkes:group_algebra}
mentioned earlier.

Cossey, Hawkes and Mann deduced from their result that $G$ is the direct product of a $p$-subgroup
and a $p'$-subgroup if and only if
$|\mathcal{S}_{p}(G)|_p=|G|_p$.
Thus the latter condition, which is the natural analogue for class sizes of the $p$-nilpotency condition
$u_{p}(G)_p=|G|_p$ of Cossey and Hawkes for character degrees~\cite[Theorem~1]{Cossey-Hawkes:group_algebra},
is equivalent to a much stronger structural property of $G$.
It is clear that adding to this the appropriate information on $|\mathcal{S}_{p}(G)|_{p'}$ yields nothing more,
namely, $G$ is the direct product of a $p$-subgroup and a $p'$-subgroup if and only if
$|\mathcal{S}_{p}(G)|=|G|_p\cdot|\mathbf{Z}(G)|_{p'}$.
This formulation can be regarded as the analogue of our Theorem~\ref{thm:direct-product-Cossey-Hawkes}.

In Section~\ref{sec:class-sizes} we provide analogues for class sizes of the results mentioned earlier which involve
character degrees prime to $p$.
Theorem~\ref{thm:Isaacs-like-class-sizes} characterizes the groups satisfying
$|\mathcal{S}_{p'}(G)|_p=|\mathbf{Z}(G)|_p$,
which is the dual version of the condition $u_{p'}(G)_p=|G:G'|_p$
appearing in Isaacs' $p$-nilpotency result.
It turns out that these are the groups $G$ where the centre of a Sylow $p$-subgroup is central in $G$.
Finally, Theorem~\ref{thm:direct-product-class-sizes} is analogous to Theorem~\ref{thm:direct-product-Isaacs}
and states that a group $G$ is the direct product of a $p$-group and a $p'$-group if and only if
$|\mathcal{S}_{p'}(G)|=|G|_{p'}\cdot|\mathbf{Z}(G)|_{p}$.

\bigskip
Part of this work (specifically, the results in Section~\ref{sec:char-degrees})
was done in Spring 1990 when I had just begun a PhD
programme at Warwick under the supervision of Trevor Hawkes.
I fondly remember those years and I am grateful to Trevor
and to the Maths Institute at Warwick for the great mathematical environment provided.

I am indebted to the referee for correcting an error in the original proof of Theorem~3.3.

%----------------------------------------------------------------
\section{Character degrees}\label{sec:char-degrees}

If $G$ is the direct product of a $\pi$-group and a $\pi'$-group, then we have
$u_{\pi'}(G)=|G|_{\pi'}\cdot|G:G'|_{\pi}$
(and, symmetrically,
$u_{\pi}(G)=|G|_{\pi}\cdot|G:G'|_{\pi'}$).
This follows at once from the fact that the irreducible characters of $G$ are exactly
all those of the form $\varphi\times\psi$,
with $\varphi$ and $\psi$ ranging over the irreducible characters of the Hall $\pi$- and, respectively,
$\pi'$-subgroup of $G$~(with notation as in \cite[Theorem~4.21]{Isa}).
At least when $\pi=\{p\}$ the converse also holds.

\begin{theorem}\label{thm:direct-product-Isaacs}
The group $G$ is the direct product of a $p$-group and a $p'$-group if and only if
$u_{p'}(G)=|G|_{p'}\cdot|G:G'|_{p}$.
\end{theorem}

\begin{proof}
Suppose that $u_{p'}(G)=|G|_{p'}\cdot|G:G'|_{p}$.
In particular,
$G$ has a normal $p$-complement $N$
according to Isaacs' result~\cite[Theorem~1]{Isaacs:group_algebra}.
Therefore, $G$ is a semidirect product $N\rtimes P$, where $P$ is any Sylow $p$-subgroup of $G$.
It remains to prove that $P$ centralizes $N$.

The proof of Isaacs' result also shows that
$u_{p'}(G)=|G:G'|_p\cdot\sum_{\theta\in\mathcal{V}}\theta(1)^2$,
where
\[
\mathcal{V}=
\{\theta\in\Irr(N)\mid
p\nmid\theta(1),\ T_G(\theta)=G\}
\]
which is just
$\{\theta\in\Irr(N)\mid
T_G(\theta)=G\}$
here, because $N$ is a $p'$-group.
Therefore, we have
$\sum_{\theta\in\mathcal{V}}\theta(1)^2=
|G|_{p'}=|N|=
\sum_{\theta\in\Irr(N)}\theta(1)^2$
and, consequently, all irreducible characters of $N$ are invariant in $G$.
Since the irreducible characters of $N$ separate the conjugacy classes of $N$ (or
according to Brauer's permutation lemma~\cite[(6.32)]{Isa}) we have
$\mathcal{K}^g=\mathcal{K}$ for every conjugacy class $\mathcal{K}$ of $N$.

A standard fact on coprime actions allows one to conclude the proof.
In fact, since $p\nmid|\mathcal{K}|$ every conjugacy class $\mathcal{K}$ of $N$
contains a $P$-orbit of length one under conjugation, and hence
$\mathcal{K}$ intersects the centralizer $\mathbf{C}_N(P)$ nontrivially.
Therefore, $N$ is the union of all conjugates of $\mathbf{C}_N(P)$ in $N$.
Since a finite group cannot be the union of all the conjugates of a proper subgroup
(see the proof of Theorem~\ref{thm:direct-product-class-sizes} for a generalization of this standard argument),
it follows that $\mathbf{C}_N(P)=N$, as desired.
\end{proof}

On the other hand, the necessary condition that
$u_{p}(G)=|G|_{p}\cdot|G:G'|_{p'}$ is not sufficient to
guarantee that $G$ is the direct product of a $p$-group and a $p'$-group.
The smallest counterexample is the holomorph of a group
of order 7, %Frobenius group $G$ of order $42$,
which has $6$ linear characters and a unique character of degree $6$, whence
$u_{2}(G)=6=|G:G'|\cdot|G'|_{2}$
(or $u_{3}(G)=6=|G:G'|\cdot|G'|_{3}$, if we prefer).
All groups satisfying that condition on $u_p(G)$ are described in the following theorem.

\begin{theorem}\label{thm:direct-product-Cossey-Hawkes}
The condition $u_{p}(G)=|G|_{p}\cdot|G:G'|_{p'}$ holds if and only if
$G$ has a normal $p$-complement $N$ with $[N,G]=N'$.
\end{theorem}

\begin{proof}
Suppose first that $u_{p}(G)=|G|_{p}\cdot|G:G'|_{p'}$.
According to Cossey and Hawkes' result~\cite[Theorem~1]{Cossey-Hawkes:group_algebra},
$G$ has a normal $p$-complement $N$.
Thus, $G=N\rtimes P$, where $P$ is a Sylow $p$-subgroup of $G$.
In this situation, the irreducible characters of $G$ with degree a power of $p$
are exactly those whose kernel contains $N'$.
In fact, if $\chi\in\Irr(G)$ and $\theta$ is an irreducible constituent of $\chi_N$,
then $\chi(1)/\theta(1)$ divides $|G:N|$, according to~\cite[(11.29)]{Isa}
(but the special case~\cite[Problem~(6.7)]{Isa} suffices here, $G/N$ being soluble).
Therefore, $\chi(1)$ is a power of $p$ if and only if $\theta(1)$ is.
However, $\theta(1)$ must also divide $|N|$, a $p'$-number.
We conclude that the irreducible characters of $G$ with degree a power of $p$
are those for which $\chi_N$ is a sum of linear characters.
Our claim follows, and hence $u_{p}(G)=|G:N'|$.
Our hypothesis implies that $|G:G'|_{p'}=|N:N'|$.
But
$G'=N'[N,P]P'$,
and hence
$|G:G'|_{p'}=|N:N'[N,P]|$.
Consequently, we have
$[N,G]=N'[N,P]=N'$,
as desired.

Now we prove the converse implication.
Thus, assume that $G$ has a normal $p$-complement $N$ and that $[N,G]=N'$.
If $P$ is any Sylow $p$-subgroup of $G$, then
$G=N\rtimes P$, and $[N,P]\le N'$.
Again according to Cossey and Hawkes' result we have $u_{p}(G)_p=|G|_p=|G:N|$.
On the other hand, we know that $u_{p}(G)=|G:N'|$, and hence
$u_{p}(G)_{p'}=|N:N'|=|N:N'[N,P]|=|G:G'|_{p'}$.
\end{proof}

\begin{cor}\label{cor:direct-product-Cossey-Hawkes}
Suppose that $u_{p}(G)=|G|_{p}\cdot|G:G'|_{p'}$ and that every $q$-element commutes with every $r$-element,
for each distinct primes $q$ and $r$ different from $p$.
Then $G$ is the direct product of a $p$-group and a $p'$-group.
\end{cor}

\begin{proof}
According to Theorem~\ref{thm:direct-product-Cossey-Hawkes}
we have $G=N\rtimes P$, where $N$ is the normal $p$-complement for $G$,
$P$ is any Sylow $p$-subgroup of $G$, and $[N,P]\le N'$.
Our assumptions imply that $N$ is nilpotent, hence its Frattini subgroup $\Phi(N)$ contains $N'$.
Each element of $P$ acts by conjugation on $N$ as an automorphism of coprime order
and induces the identity automorphism on $N/\Phi(N)$.
Therefore, it acts trivially on $N$ according
to~\cite[Satz~3.18]{Hup}.
We conclude that $G$ is the direct product of $N$ and $P$.
\end{proof}

\begin{rem}\label{rem:char-degrees}
Suppose that $|G|$ and $|G:G'|$ are given.
Isaacs' criterion~\cite[Theorem~1]{Isaacs:group_algebra}
asserts that $G$ is $p$-nilpotent exactly if $u_{p'}(G)_p$ is as {\em small} as possible,
while Cossey and Hawkes' criterion~\cite[Theorem~1]{Cossey-Hawkes:group_algebra}
says that $G$ is $p$-nilpotent exactly if $u_{p}(G)_p$ is as {\em large} as possible.

The proof of Theorem~\ref{thm:direct-product-Isaacs} shows that, in general,
if $u_{p'}(G)_{p}=|G:G'|_{p}$ then
$u_{p'}(G)_{p'}\le|G|_{p'}$.
Hence $G$ is the direct product of a $p$-group and a $p'$-group exactly if
$u_{p'}(G)_p$ is as small as possible and
$u_{p'}(G)_{p'}$ is as large as possible subject to the former condition.
Similarly, the proof of Theorem~\ref{thm:direct-product-Cossey-Hawkes} shows that, in general,
if $u_{p}(G)_{p}=|G|_{p}$ then
$u_{p}(G)_{p'}\ge|G:G'|_{p'}$
(and actually more, namely $|G:G'|_{p'}$ divides $u_{p}(G)_{p'}$).
Therefore, $G$ has a normal $p$-complement $N$ with $[N,G]=N'$ exactly if
$u_{p}(G)_p$ is as large as possible and
$u_{p}(G)_{p'}$ is as small as possible subject to the former condition.
\end{rem}

%----------------------------------------------------------------
\section{Class sizes}\label{sec:class-sizes}

The following theorem is an analogue for class sizes of Isaacs' characterization
of $p$-nilpotency in terms of $u_{p'}(G)$~\cite[Theorem~1]{Isaacs:group_algebra}.
The proof also resembles Isaacs' proof rather closely.
Curiously, however, the group-theoretic equivalent condition in Theorem~\ref{thm:Isaacs-like-class-sizes} looks somehow dual
to that in Theorem~\ref{thm:direct-product-Cossey-Hawkes},
which is based on Cossey and Hawkes' result and involves $u_p(G)$.

\begin{theorem}\label{thm:Isaacs-like-class-sizes}
For every group $G$ the order of its centre $\mathbf{Z}(G)$ divides $|\mathcal{S}_{p'}(G)|$.
Furthermore, $|\mathcal{S}_{p'}(G)|_p=|\mathbf{Z}(G)|_p$ if and only if
$[\mathbf{Z}(P),G]=1$, where $P$ is a Sylow $p$-subgroup of $G$.
\end{theorem}

\begin{proof}
An element of $G$ has class size a $p'$-number if and only if its centralizer contains a Sylow $p$-subgroup of $G$.
Therefore, we have
$\mathcal{S}_{p'}(G)=\bigcup_{P\in\mathrm{Syl}_p(G)}\mathbf{C}_G(P)$.
Since $\mathbf{Z}(G)$ is contained in $\mathbf{C}_G(P)$ for all $P\in\mathrm{Syl}_p(G)$
it follows that $\mathcal{S}_{p'}(G)$ is a union of cosets of $\mathbf{Z}(G)$,
and hence $|\mathbf{Z}(G)|$ divides $|\mathcal{S}_{p'}(G)|$.
Actually, one can prove a stronger conclusion, which is not needed in the present proof, noting that
$\mathbf{C}_G(\mathbf{O}^{p'}(G))=\bigcap_{P\in\mathrm{Syl}_p(G)}\mathbf{C}_G(P)$
since
$\mathbf{O}^{p'}(G)=\langle P: P\in\mathrm{Syl}_p(G)\rangle$.
Thus, $\mathcal{S}_{p'}(G)$ is a union of cosets of $\mathbf{C}_G(\mathbf{O}^{p'}(G))$,
and so $|\mathbf{C}_G(\mathbf{O}^{p'}(G))|$ divides $|\mathcal{S}_{p'}(G)|$.

Fix a Sylow $p$-subgroup $P$ of $G$ and let $Z=P\cap\mathbf{Z}(G)$, whence $Z\in\mathrm{Syl}_p(\mathbf{Z}(G))$.
Let $\overline{\mathcal{S}}$ be the image of $\mathcal{S}_{p'}(G)$ in $G/Z$.
Note in passing that $\overline{\mathcal{S}}$ is contained in $\mathcal{S}_{p'}(G/Z)$
because $\mathbf{C}_{G/Z}(gZ)\ge\mathbf{C}_{G}(g)/Z$ for every $g\in G$,
but $\overline{\mathcal{S}}$ may well be a proper subset of $\mathcal{S}_{p'}(G/Z)$.
However, we claim that for $g\in\mathcal{S}_{p'}(G)$ we have $\mathbf{C}_{G/Z}(gZ)=\mathbf{C}_{G}(g)/Z$,
which is the same as $\mathbf{N}_G(gZ)=\mathbf{C}_{G}(g)$.
In fact, because $Z$ is a central subgroup of $G$, the map $\mathbf{N}_G(gZ)\to Z$ given by $x\mapsto[g,x]$
is a group homomorphism with kernel $\mathbf{C}_{G}(g)$.
Since $Z$ is a $p$-group this can only be the trivial homomorphism
when $|G:\mathbf{C}_{G}(g)|$ is a $p'$-number, that is, when $g\in\mathcal{S}_{p'}(G)$.
This proves our claim, which is also equivalent to saying that
the elements of $gZ$ are pairwise nonconjugate in $G$.

Consider the action of the $p$-group $P/Z$ on $\overline{\mathcal{S}}$ by conjugation.
The fixed points of this action are exactly the elements of $\mathbf{C}_G(P)/Z$.
In fact, if $gZ$ is centralized by $P/Z$ then $g$ is centralized by $P$
because of what we have proved in the previous paragraph.
The converse is obvious.
It follows that
\[
|\mathcal{S}_{p'}(G)|/|\mathbf{Z}(G)|_p=
|\overline{\mathcal{S}}|\equiv
|\mathbf{C}_{G}(P)/Z|\pmod{p}.
\]
Consequently, $|\mathcal{S}_{p'}(G)|_p=|\mathbf{Z}(G)|_p$
holds if and only if
$\mathbf{C}_{G}(P)/Z$ is a $p'$-group.
This occurs if and only if
$\mathbf{Z}(P)\le\mathbf{Z}(G)$,
because $\mathbf{Z}(P)$ is the Sylow $p$-subgroup of $\mathbf{C}_{G}(P)$.
\end{proof}

Another way of stating the condition $[\mathbf{Z}(P),G]=1$ in Theorem~\ref{thm:Isaacs-like-class-sizes}
is that the only $p$-elements with class length a $p'$-number are central.

\begin{rem}
The first statement of Theorem~\ref{thm:Isaacs-like-class-sizes} generalizes to a set of primes $\pi$
instead of $p$.
In fact, since
$\mathcal{S}_{\pi'}(G)=\bigcap_{p\in\pi}\bigcup_{P\in\mathrm{Syl}_p(G)}\mathbf{C}_G(P)$, we have that
$\mathcal{S}_{\pi'}(G)$ is a union of cosets of $\mathbf{Z}(G)$, and so
$|\mathbf{Z}(G)|$ divides $|\mathcal{S}_{\pi'}(G)|$.
Similarly, Isaacs' proof of~\cite[Theorem~1]{Isaacs:group_algebra} shows, more generally, that
$|G:G'|_\pi$ divides $u_{\pi'}(G)$ for every finite group $G$ and set of primes $\pi$.
However, the second parts of both Theorem~\ref{thm:Isaacs-like-class-sizes} and
Isaacs' Theorem~\cite[Theorem~1]{Isaacs:group_algebra} do not admit obvious
generalizations to a set of primes $\pi$.
\end{rem}

It is easy to verify that if $G$ is the direct product of a $\pi$-group and a $\pi'$-group, then we have
$|\mathcal{S}_{\pi'}(G)|=|G|_{\pi'}\cdot|\mathbf{Z}(G)|_{\pi}$.
In complete analogy with the case of character degrees examined in Section~\ref{sec:char-degrees},
when $\pi=\{p\}$ the converse holds.
The proof is very similar to the proof of its dual result, Theorem~\ref{thm:direct-product-Isaacs}.

\begin{theorem}\label{thm:direct-product-class-sizes}
The group $G$ is the direct product of a $p$-group and a $p'$-group if and only if
$|\mathcal{S}_{p'}(G)|=|G|_{p'}\cdot|\mathbf{Z}(G)|_{p}$.
\end{theorem}

\begin{proof}
Assume that $|\mathcal{S}_{p'}(G)|=|G|_{p'}\cdot|\mathbf{Z}(G)|_{p}$, let $P$
be a Sylow $p$-subgroup of $G$, and set $Z=\mathbf{Z}(P)$.
According to Theorem~\ref{thm:Isaacs-like-class-sizes} we have
$Z\le\mathbf{Z}(G)$.
We have
$\mathcal{S}_{p'}(G)=\bigcup_{x}\mathbf{C}_G(P)^x$,
with $x$ ranging over a right transversal for
$\mathbf{N}_G(\mathbf{C}_G(P))$ in $G$.
Since $P\le\mathbf{N}_G(\mathbf{C}_G(P))$,
the cardinality $|G:\mathbf{N}_G(\mathbf{C}_G(P))|$ of such a transversal is a $p'$-number.
Furthermore, $|\mathbf{C}_{G}(P):Z|$ is a $p'$-number,
because $Z$ is the Sylow $p$-subgroup of $\mathbf{C}_{G}(P)$.

Now we apply a generalization of the standard argument used to prove that a finite group cannot be the union of all
the conjugates of a proper subgroup.
We have
\begin{align*}
|G|_{p'}\cdot|Z|=|\mathcal{S}_{p'}(G)|
&\le
|Z|+(|\mathbf{C}_G(P)|-|Z|)\cdot|G:\mathbf{N}_G(\mathbf{C}_G(P))|
\\&\le
|Z|+|G|_{p'}\cdot|Z|-|G:\mathbf{N}_G(\mathbf{C}_G(P))|\cdot|Z|,
\end{align*}
because the product of $|\mathbf{C}_{G}(P):Z|$ and $|G:\mathbf{N}_G(\mathbf{C}_G(P))|$ is a $p'$-number
dividing $|G|$.
We deduce that $|G:\mathbf{N}_G(\mathbf{C}_G(P))|=1$, that is, $\mathbf{C}_G(P)$ is normal in $G$.
Consequently, $\mathbf{C}_G(P)$ equals $\mathcal{S}_{p'}(G)$, and then by hypothesis,
its index in $G$ equals $|G:\mathbf{Z}(G)|_p$, a power of $p$.
It follows that $P\mathbf{C}_G(P)=G$ and, therefore, $P$ is a normal subgroup of $G$.

According to the Schur-Zassenhaus theorem $P$ has a complement $H$ in $G$,
and hence $G=P\rtimes H$.
Since $\mathcal{S}_{p'}(G)$ equals $\mathbf{C}_G(P)$,
our hypothesis implies that $|G:\mathbf{C}_G(P)|$ is a power of $p$.
Therefore $\mathbf{C}_G(P)$, being a normal subgroup of $G$, contains every Sylow $q$-subgroup of $G$ for $q\neq p$,
and hence contains $H$.
We conclude that $G$ is the direct product of $P$ and $H$.
\end{proof}

\begin{rem}\label{rem:class-lengths}
We attempt a comparison with Remark~\ref{rem:char-degrees}, assuming that $|G|$ and $|\mathbf{Z}(G)|$ are given.
Theorem~\ref{thm:Isaacs-like-class-sizes} says that
$\mathbf{Z}(P)\le\mathbf{Z}(G)$ for a Sylow $p$-subgroup $P$ of $G$
exactly if $|\mathcal{S}_{p'}(G)|_p$ is as small as possible.
On the other hand, the result of Cossey, Hawkes and Mann~\cite{Cossey-Hawkes-Mann}
mentioned in the Introduction implies that
$G$ is the direct product of a $p$-group and a $p'$-group
exactly if
$|\mathcal{S}_{p}(G)|_p$ is as large as possible.

So far the similarity with the case of character degrees is not striking, but becomes more so
if we include also the $p'$-parts of our invariants into consideration.
The proof of Theorem~\ref{thm:direct-product-class-sizes} shows that, if
$|\mathcal{S}_{p'}(G)|_{p}=|\mathbf{Z}(G)|_{p}$ then
$|\mathcal{S}_{p'}(G)|\le|G|_{p'}\cdot|\mathbf{Z}(G)|_{p}$ in general, and
so $|\mathcal{S}_{p'}(G)|_{p'}\le|G|_{p'}$.
Hence $G$ is the direct product of a $p$-group and a $p'$-group exactly if
$|\mathcal{S}_{p'}(G)|_p$ is as small as possible and $|\mathcal{S}_{p'}(G)|_{p'}$
is as large as possible subject to the former
condition.
Symmetrically, it is correct to say that $G$ is the direct product of a $p$-group and a $p'$-group exactly if
$|\mathcal{S}_{p}(G)|_p$ is as large as possible and
$|\mathcal{S}_{p}(G)|_{p'}$ is as small as possible subject to the former condition.
Of course, the latter condition is redundant here,
since the former condition already forces $G$ to be the direct product of a $p$-group and a $p'$-group,
whence $|\mathcal{S}_{p}(G)|_{p'}=|\mathbf{Z}(G)|_{p'}$.
\end{rem}

\bibliography{References}

\end{document}